\documentclass[11pt]{article}
\pdfoutput=1
\usepackage{amsmath}
\usepackage[pdftex]{graphicx}
\usepackage{caption}
\usepackage{float}
\usepackage{hyperref}
\usepackage{longtable}
\usepackage{tabularx}
\usepackage{booktabs}
\usepackage{authblk}
\usepackage[margin=3.5cm]{geometry}

\restylefloat{figure}

\newcommand{\head}[1]{\textnormal{\textbf{#1}}}

\usepackage[sort]{natbib}

\title{Multi-Haul Quasi Network Flow Model for \\ Vertical Alignment Optimization
\footnote{The Version  of Record of this manuscript has been published and is available in Engineering Optimization (GENO); Estimated publication date: Jan 12, 2017 (online); \url{http://dx.doi.org/10.1080/0305215X.2016.1271880}}
}
\author[1]{Vahid Beiranvand}
\author[2]{Warren Hare}
\author[1]{Yves Lucet\footnote{corresponding author: yves.lucet@ubc.ca}}
\author[1]{Shahadat Hossain}
\affil[1]{\small Computer Science, University of British Columbia, Kelowna, BC, Canada}
\affil[2]{Mathematics, University of British Columbia, Kelowna, BC, Canada}
\date{\today}

\begin {document}
\maketitle
\begin{abstract}
\noindent The vertical alignment optimization problem for road design aims to generate a vertical alignment of a new road with a minimum cost, while satisfying safety and design constraints. We present a new model called multi-haul quasi network flow (MH-QNF) for vertical alignment optimization that improves the accuracy and reliability of previous mixed integer linear programming models. We evaluate the performance of the new model compared to two state-of-the-art models in the field: the complete transportation graph (CTG) and the quasi network flow (QNF) models. The numerical results show that, within a 1\% relative error, the proposed model is robust and solves more than 93\% of test problems compared to 82\% for the CTG and none for the QNF. Moreover, the MH-QNF model solves the problems approximately 8 times faster than the CTG model.
\end{abstract}

\section{Introduction}

Transportation infrastructures are an important sign of development and welfare in a given country \cite{38, 39}. While planners normally try to balance and coordinate among various transportation methods, roads and highways are considered the leading component of transportation infrastructures. Due to the fast increase of traffic volume and loadings, high costs of road construction projects, and environmental and safety impacts, road construction require continuing innovation. Road design is one of the early and significant steps in road construction projects and it influences the construction cost and other contributing constraints significantly.

Road design starts by selecting the corridor within which the road is to be constructed.  From here, a so called horizontal alignment is fixed.  Next vertical alignment and earthwork movement problems must be solved. In the literature, these last two steps are often combined for the sake of efficiency \cite{6, 22}.  Some studies consider the simultaneous optimization of the horizontal alignment, vertical alignment, and earthwork movement~\cite{12, 44, c6, 45}. Our motivation to focus solely on the last two problems is threefold. First, the horizontal alignment problem is very complex and can currently only be solved by heuristics that do not guarantee global optimality. Moreover, the resulting roads generally still require a final detailed optimization of the vertical alignment and earthwork movement problems~\cite{LEE-09}. Next, when upgrading a road, the horizontal alignment is already decided and only the vertical and earthwork problems need to be solved. Finally, building a deterministic global optimization algorithm to solve the combined vertical alignment and earthwork problem is a serious step towards solving the complete problem.  Specifically, such algorithms can be used in bilevel formulations that are designed to simultaneously solve horizontal, vertical and earthwork problems, as proposed in~\cite{MONDAL-15}.

Heuristic optimization methods, such as Genetic Algorithms (GAs), have been used to solve the road alignment problem. \Citet{18} developed a mathematical model for vertical alignment and proposed a heuristic method to solve it. \Citet{37} present a comprehensive review of the road design problem and the related challenges. They also present solution algorithms based on GAs to optimize highway alignments by adapting operators and encoding schemes. Other heuristic-based optimization methods for highway alignment optimization include \cite{40, 41, 42, 43, 8}.

Another early and popular approach relied on dynamic programming. Dynamic programming was first employed \cite{21} to develop two models for vertical alignment. They use piecewise-linear segments to represent the alignment. Based on the results of the paper, the state parametrization model is more flexible to be used for solving a three-dimensional problem while the dynamic programming model has a better performance in terms of computational time and ease of formulation. In general, compared to optimizing the horizontal alignment, a dynamic programming approach has a better performance in vertical alignment optimization; however, the resulting alignment is limited to a finite set of points at each station. Hence, only a portion of the problem search space is explored \cite{28}. Other models using dynamic programming include~\cite{29, 30, 46}.

Our focus is to minimize the road construction costs (primarily excavation, embankment, and hauling costs), subject to road safety and quality constraints. While costs of land, pavement, vehicle operation, motorist time, accidents, and environmental impacts all contribute to the cost, especially when one considers highways, they have less impact when the horizontal alignment is fixed, or for logging roads. Optimizing solely road construction costs allows us to simplify our model and to take advantage of specific problem structures to speed up the computation.

The current state-of-the-art on the vertical alignment and earthwork movement problem relies on linear programming, or mixed integer linear programming.  \Citet{14} proposed a two-step linear programming model to optimize vertical alignment and earthwork movement.  In the first step by enumerating all feasible grades, a feasible road grade is chosen and for each grade the amount of cut and fill is calculated. In the second step the earthwork transportation cost is minimized using linear programming. This trial and error approach does not guarantee finding the global solution. \Citet{15} proposed a linear programming model for vertical alignment which incorporates the grade selection and earthwork allocation stages. \Citet{22} extended the work of \citet{14} to handle the sharp connectivity points problem in piecewise linear models by representing the road profile as a continuous quadratic curve. \Citet{33} proposed a model in which the road profile is represented as a one-dimensional spline. This model improves the computational efficiency and guarantee the global optimality by solving a single linear programming problem. \Citet{13} improved the \cite{33} model by introducing new gap and slope constraints that reduce errors in the model.

\Citet{6} developed a mixed integer linear program (MILP) model for vertical alignment in which the blocks are also taken into account. In \cite{5} road side slopes are added to the previous MILP model proposed in \cite{6} which yields a more accurate earthwork cost calculation. \Citet{4} also improved the performance of the MILP model (called complete transportation graph or CTG) and studied a novel quasi network flow (QNF) model for vertical alignment optimization considering the earthwork allocations. The QNF model presented significant improvements in computational time, but at a cost of decreased accuracy in the final solution.

In this paper, we improve the QNF model of \cite{4}, by presenting a new model called the multi-haul cost quasi network flow (MH-QNF) model.  The MH-QNF model balances the accuracy of the CTG model in approximating the earthwork cost with the low computational cost of the QNF model. The main idea behind the new model is to differentiate between hauls of different distances in earthwork allocation. For simplicity, we consider three possible distances, short hauls, middle hauls, and long hauls; however, any number of haul distances can be accommodated by the model.

\section{Model description}

The road is approximated as a quadratic spline. It is split into $m$ segments indexed by $\mathcal{G}=\{1,2,...,m\}$. For all $g \in \mathcal{G}$, each segment is represented using the following equation

\begin{equation}
P_g(s)=a_{g,1}+a_{g,2}s+a_{g,3}s^2
\end{equation}
in which $s$ is the distance from the beginning of the road.

Each segment is further split into sections so that the $g$\textsuperscript{th} spline segment is made of $n_g$ sections indexed by the set $\mathcal{S}_g=\{ 1,2,...,n_g \}$. The total number of sections in a road is $n=\sum_{g \in \mathcal{G}}n_g$ and these sections are indexed by the set $\mathcal{S}=\{ 1,2,...,n \}$. The section indexes of a specific spline segment $g$ are mapped to the actual section index set $\mathcal{S}$ using the function $\varphi : (\mathcal{G}\times\mathcal{S}_g) \rightarrow \mathcal{S}$. For example, if $\varphi(g,j) = i $ then  $s_i=s_{\varphi(g,j)}$ for all $i\in \mathcal{S}, g \in \mathcal{G}, j \in \mathcal{S}_g$.
So the spline function representing the road profile is
\[P(s) =
	\begin{cases}
      P_1(s) & \textbf{\text{if }} s_{\varphi(1,1)} \le s \le s_{\varphi(1,n_1)}, \\
      P_2(s) & \textbf{\text{if }} s_{\varphi(2,1)} \le s \le s_{\varphi(2,n_2)}, \\
      \vdots \\
      P_{n_g}(s) & \textbf{\text{if }} s_{\varphi(n_g,1)} \le s \le s_{\varphi(2,n_{n_g})}.
	\end{cases}
\]

In the vertical alignment problem, one goal is to obtain the optimal offset between the ground profile and the road profile. In addition, the required cut and fill volumes should be calculated for each offset. The offset is denoted by $u_{i}, i \in \mathcal{S}$, and the cut and fill volumes of a section $i$ are denoted by $V^{+}_{i}$ and $V^{-}_{i}, i \in \mathcal{S}$ respectively. (Shrinkage and swell factors for each material can be easily incorporated in the model but are not included to simplify the presentation.)

When constructing a road we need to consider borrow pits to bring in material and waste pits to dump extra materials. They are modeled as external sections and indexed by the sets $\mathcal{B}=\{1, 2,...,n_\beta\}$ and $\mathcal{W}=\{ 1, 2,...,n_w \}$ in which $n_\beta$ is the number of borrow pit sections and $n_w$ is the number of waste pit sections. A borrow pit index is mapped to the section index to which it is attached by the function $\vartheta:\mathcal{B} \rightarrow \mathcal{S}$. Similarly, the function $\delta:\mathcal{W} \rightarrow \mathcal{S}$ maps the waste pit index to the section index to which it is attached.

Access roads are required in road construction and they are used as gateways to the road being constructed. Access roads are linked to a section of the road and are indexed by the set $\mathcal{R}=\{ 1, 2,...,n_r \}$ in which $n_r$ is the number of access roads in the construction site. The function $\varrho:\mathcal{R} \rightarrow \mathcal{S}$ maps an access road index to the section index to which it is linked. The capacity of the $i$\textsuperscript{th} borrow pit (respectively waste pit) is denoted by $C^b_i$ (respectively $C^w_i$). The set $\mathcal{N} = \mathcal{S} \cup \mathcal{B} \cup \mathcal{W}=\{ 1,2,...,n+n_\beta+n_w \}$ represents all the indexes for sections, borrow pits, and waste pits. The {\it dead haul distance} is the distance between a borrow/waste pit and the section to which it is linked. The dead haul distance of the $i$\textsuperscript{th} borrow or waste pit is denoted by $\tilde{d}_i$.

The side slope is defined as the steady decrease/increase of height when moving orthogonally from the road profile to the ground profile in a cut/fill section. Side slopes are represented as trapezoid shaped cross-sections and are approximated using stacked rectangles to preserve the linearity of the model, see \cite{4}. The model can handle nonsymmetric section-dependent side slopes to account for local soil type and sloped terrains as well as cuts and fills at the same section.

Another important consideration in vertical alignment optimization are blocks \cite{6}. Blocks are obstacles that need to be removed to access some parts of the corridor. They indicate a river or mountain, and require building a bridge or a tunnel. We define $n_b$ as the number of blocks in the corridor. Blocks are indexed by the set $\mathcal{I}=\{ 1, 2,...,n_b \}$. The function $\gamma:\mathcal{I} \rightarrow \mathcal{S}$ maps a block index to the section index to which it is linked. To model the block removal process we use a time step $t$ which specifies the time at which a block is removed. In \cite{6}, it is shown that to remove $n_b$ blocks we need at most $n_b+1$ time steps. So we define the set $\mathcal{T}=\{ 0, 1, 2,...,n_b \}$ to model the required time steps and we use the binary variables $y_{kt}$ for each block $k \in \mathcal{I}$ and time step $t \in \mathcal{T}$ to determine whether a block is removed or not. After a block is removed we can move material over its section, see constraints \eqref{14}--\eqref{21}.

In the MH-QNF model, similarly to the QNF model \cite{4}, we draw sections and pits as nodes while arcs show feasible movements between them. In the QNF model the authors use a single haul to move materials. This is equivalent to assuming a single type of earthmoving equipment.  In real roadway construction sites, different equipments are required for earthmoving tasks. For example, for the short distances using a bulldozer may be more economical, while for a long distance movement a truck may be preferred. Therefore, to obtain a more realistic model, we extend the previous QNF model \cite{4} by using multiple hauling paths, which is more realistic and yields a more precise earthwork solution.

\begin{figure}[ht]
    \centering
	\includegraphics[scale=0.7]{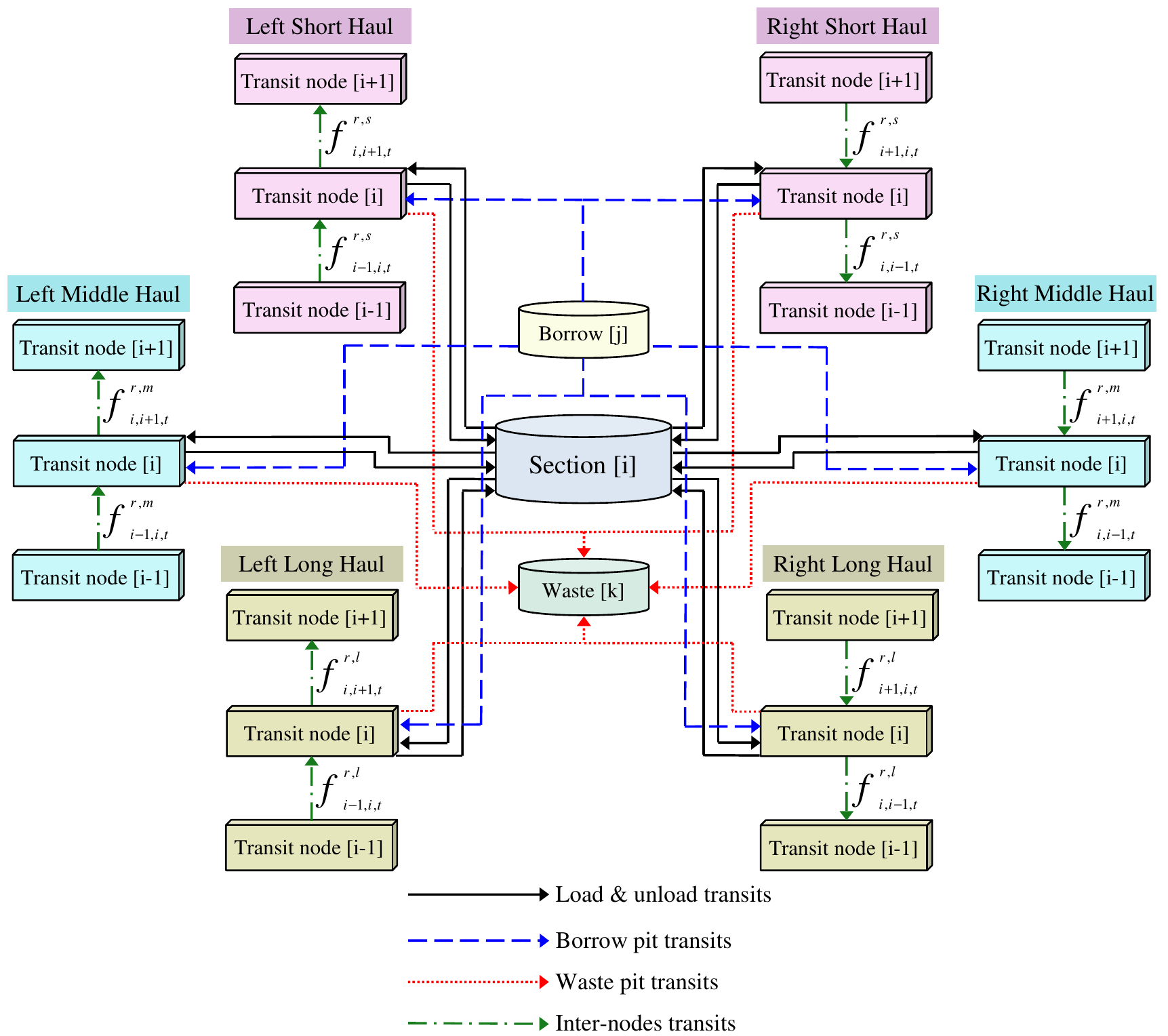}
    \caption{The general scheme of a typical section in the MH-QNF model}
    \label{f:1}
\end{figure}

Figure \ref{f:1} shows the typical flows for a section in the proposed MH-QNF model.

In this model, the material can be moved from the current section $i$ to the next section ($i+1$) or the previous one ($i-1$). Depending on the distance of a hauling task, the material movement can be performed through one of a number of hauling paths indexed by $\mathcal{H}=\{ 1, 2, ..., n_h\}$.  (In our figures and numerical tests we shall apply $n_h=3$, which is based on consultation with our industry partner, whose standard software includes {\em free-haul}, {\em over-haul}, and {\em end-haul}.)

When the flow of material reaches a section, materials can be unloaded to fill the section, a cut from the section can be performed with the new materials being added to the flow, or the section remains unchanged and the material is transferred to the next section. To transfer material to the left or right we use virtual transit nodes for both left and right directions. Since we have $n_h$ different hauls, $n_h$ groups of transit nodes are denoted in the model, see Figure \ref{f:1} (which represents three hauls for consistency with numerical experiments). If material is moved to the next section, the left-side transit nodes are employed; otherwise if the material is moved to the previous section the right-hand side transit nodes are used. In Figure \ref{f:1}, we define $f_{i,i+1,t}^{r,s}$ and $f_{i,i-1,t}^{r,s}$ (for all $i \in \mathcal{S}, t \in \mathcal{T}$) as the flow of material from the left (right) short haul transit node $i$ to the left (right) short haul transit node $i+1$ (or $i-1$) at time step $t$, respectively. Similarly, we use the variables $f_{i,i+1,t}^{r,m}$ and $f_{i,i-1,t}^{r,m}$ (for all $i \in \mathcal{S}, t \in \mathcal{T}$) for the middle hauling path, and the variables $f_{i,i+1,t}^{r,l}$ and $f_{i,i-1,t}^{r,l}$ (for all $i \in \mathcal{S}, t \in \mathcal{T}$) for the long hauling path.

\begin{figure}[ht]
    \centering
	\includegraphics[scale=0.7]{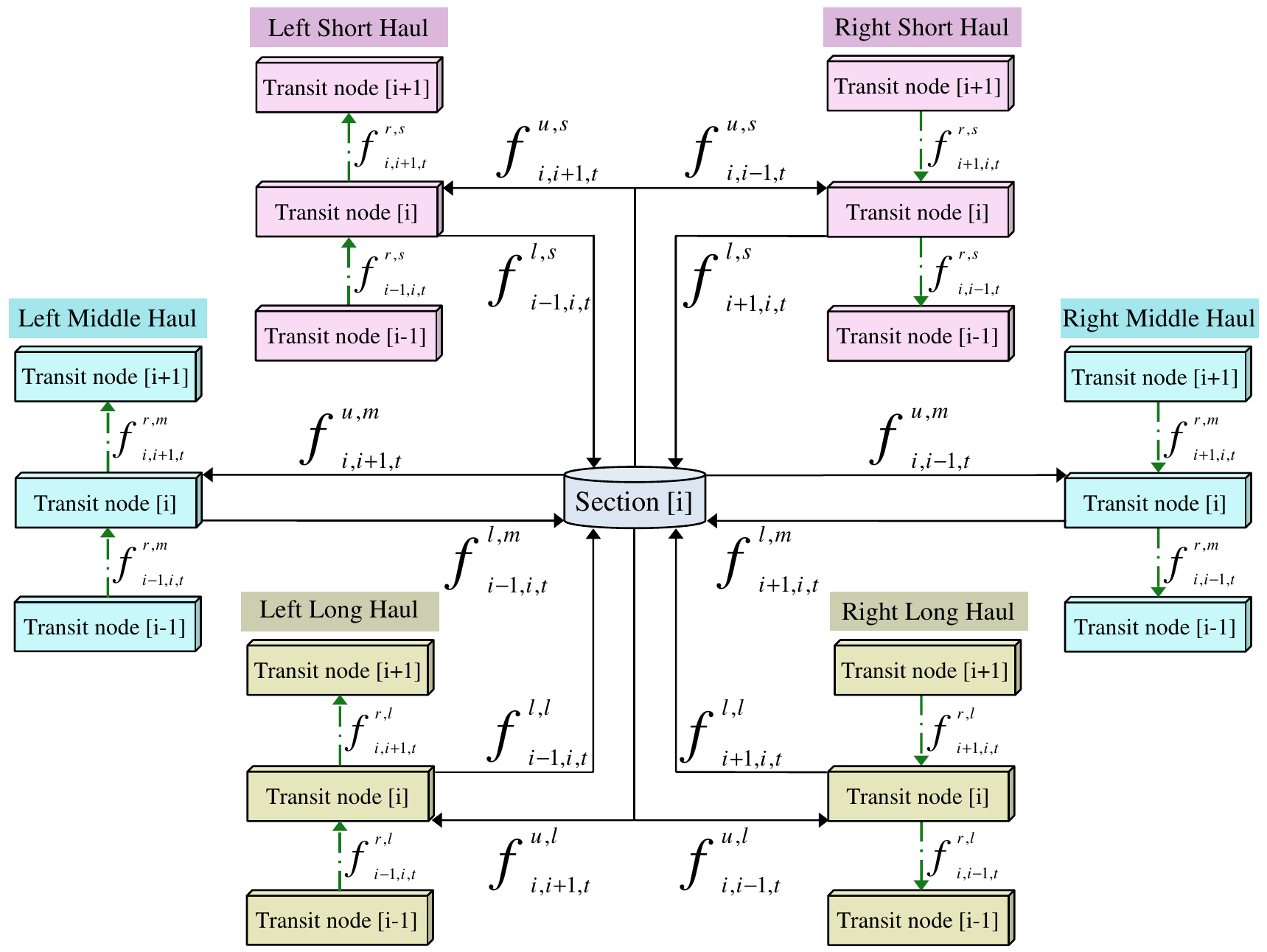}
    \caption{Load and unload flows in the MH-QNF model}
    \label{2}
\end{figure}

Figure \ref{2} shows the load and unload flows that model the fill and cut volumes of materials for each section.

After cutting the materials we can transfer them to the left or right using any of the three hauling options (short, middle, or long). Similarly, a section can be filled using the materials taken from left or right transit nodes. To model these possibilities, we introduce the variables $f_{i-1,i,t}^{l,s}$ and $f_{i+1,i,t}^{l,s}$ (for all $i \in \mathcal{S}, t\in\mathcal{T}$) as the load (fill) flows of materials from the left and right transit nodes for the short hauling path. The variables $f_{i-1,i,t}^{l,m}$, $f_{i+1,i,t}^{l,m}$; and $f_{i-1,i,t}^{l,l}$ and $f_{i+1,i,t}^{l,l}$ (for all $i \in \mathcal{S}, t\in\mathcal{T}$) are used as the load flows of materials for the middle and long hauling paths, respectively. Similarly, the variables $f_{i,i+1,t}^{u,s}$ and $f_{i,i-1,t}^{u,s}$ (for all $i \in \mathcal{S}, t\in\mathcal{T}$) are the unload (cut) flows of materials to the left and right transit nodes for the short hauling path. In the same way, the variables $f_{i,i+1,t}^{u,m}$, $f_{i,i-1,t}^{u,m}$; and $f_{i,i+1,t}^{u,l}$, $f_{i,i-1,t}^{u,l}$ (for all $i \in \mathcal{S}, t\in\mathcal{T}$) are the unload flows for the middle and long hauling paths, respectively.

Each borrow pit is attached to a section. A borrow flow denotes material movement from the pit to that specific section; from there the material may move to a transit node. The variables $f_{j,\vartheta (j)+1,t}^{b,s}$ and $f_{j,\vartheta (j)-1,t}^{b,s}$ (for all $j \in \mathcal{B}, t\in\mathcal{T}$) can be used to transfer materials from a borrow pit to the next or the previous section via the short haul path. Similarly, to represent the borrow flows for the middle and long hauling paths we define the variables $f_{j,\vartheta (j)+1,t}^{b,m}$, $f_{j,\vartheta (j)-1,t}^{b,m}$; and $f_{j,\vartheta (j)+1,t}^{b,l}$ and $f_{j,\vartheta (j)-1,t}^{b,l}$ (for all $j \in \mathcal{B}, t\in\mathcal{T}$), respectively. Waste flows are handled similarly.

A typical real-world road construction project involves several types of materials with different excavation, embankment, and hauling costs. The proposed model is able to handle more than one material. We define the index set $\mathcal{M}=\{ 1, 2, ..., n_m \}$ corresponding to the indexes of four types of materials. In the present research we chose $n_m = 4$ based on consultation with our industry partner.

The cost components of a vertical alignment problem consist of the embankment, excavation, and hauling costs of materials in the construction site. The per unit of volume excavation (cut or unload from a section) cost of material is represented by $p_1$, $p_2$, ... $p_{n_m}$ for different types of materials, respectively. Similarly, $q_1$, $q_2$,  ... $q_{n_m}$  denote the per unit volume of embankment (fill or load to a section) cost of material for the different types of materials. Finally, for the short hauling path, the hauling cost of materials from the section $i$ to the section $i-1$ (respectively $i+1$) is defined as $c_{i,i-1}^{r,s}=c_sd_{i,i-1}$ (respectively, $c_{i,i+1}^{r,s}=c_sd_{i,i+1}$) in which $c_s$ is the cost of moving one unit volume of materials per unit distance for the short hauling path, and $d_{i,j}$ is the distance between sections $i$ and $j$. Similarly, for the middle hauling path we define the hauling cost components $c_{i,i-1}^{r,m}=c_md_{i,i-1}$ and $c_{i,i+1}^{r,m}=c_md_{i,i+1}$; and for the long hauling path we have $c_{i,i-1}^{r,l}=c_ld_{i,i-1}$ and $c_{i,i+1}^{r,l}=c_ld_{i,i+1}$. There is another cost when cutting materials from a section or a borrow pit and it is the loading cost, which is equipment dependent. Therefore, we define the loading cost components $y_s$, $y_m$, and $y_l$ as the cost of loading different types of construction equipments which use short, middle, and long hauling paths, respectively.

In the MH-QNF model the objective is to minimize the total cost of excavation, embankment, and hauling tasks. Therefore, the objective function is defined as follows.

\begin{multline*}
\min
\Bigg[ \sum_{\substack{i \in \mathcal{S \cup B} \\ m \in \mathcal{M} \\ h \in \mathcal{H} }} (p_{m}+y_{h})V_i^+
       + \sum_{\substack{i \in \mathcal{S \cup W} \\ m \in \mathcal{M} }} q_mV_i^-
 + \sum_{\substack{i \in \mathcal{S} \\ h \in \mathcal{H} \\ t \in \mathcal{T} }} (c_{i,i-1}^{r,h} f_{i,i-1,t}^{r,h} + c_{i,i+1}^{r,h} f_{i,i+1,t}^{r,h} ) \\
 + \sum_{\substack{j \in \mathcal{B} \\ m \in \mathcal{M} \\ h \in \mathcal{H} \\ t \in \mathcal{T} }} (p_{m}+y_{h}+ c_h \tilde{d}_j ) (f_{j,\vartheta (j)-1,t}^{b,h}+f_{j,\vartheta (j)+1,t}^{b,h}) \\
 + \sum_{\substack{k \in \mathcal{W} \\ m \in \mathcal{M} \\ h \in \mathcal{H} \\ t \in \mathcal{T} }} (q_{m}+ c_h \tilde{d}_k) (f_{\delta (k)-1,k,t}^{w,h} + f_{\delta (k)+1,k,t}^{w,h})\Bigg].
\end{multline*}

Since the transit nodes are only intermediate virtual points that provide transits for the flows of materials, the sum of the input flows to a transit node must be equal to the sum of the output flows from the node. Therefore, for all $i\in \mathcal{S}, h \in \mathcal{H}, t \in \mathcal{T}$ the \textit{flow constraints} are defined as
\begin{align*}
    f_{i-1,i,t}^{r,h} + f_{i,i+1,t}^{u,h} + \sum_{\substack{ j \in \mathcal{B} \\ \vartheta(j)=i }} f_{j,i+1,t}^{b,h} &= f_{i,i+1,t}^{r,h} + f_{i-1,i,t}^{l,h} + \sum_{\substack{ k \in \mathcal{W} \\ \delta(k)=i }} f_{i-1,k,t}^{w,h},     \\
    f_{i+1,i,t}^{r,h} + f_{i,i-1,t}^{u,h} + \sum_{\substack{ j \in \mathcal{B} \\ \vartheta(j)=i }} f_{j,i-1,t}^{b,h} &= f_{i,i-1,t}^{r,h} + f_{i+1,i,t}^{l,h} + \sum_{\substack{ k \in \mathcal{W} \\ \delta(k)=i }} f_{i+1,k,t}^{w,h}.
\end{align*}

The section nodes are treated as the source and destination of the flows. Therefore, the total sum of load flows to a section or a waste pit should be equal to the total fill volume of materials from that section or waste pit. Similarly, the total sum of unload flows from a section or a borrow pit should be equal to the total cut volume of materials from that section or borrow pit. These constraints are called \textit{balance constraints} and defined as
\begin{align*} 
\sum_{t \in \mathcal{T}} (f_{i,i+1,t}^{u,h}+f_{i,i-1,t}^{u,h}) &= V_i^+, \quad &i \in \mathcal{S}, h \in \mathcal{H}, & \\
\sum_{\substack{ t \in \mathcal{T} }} (f_{j,\vartheta(j)+1,t}^{b,h} + f_{j,\vartheta(j)-1,t}^{b,h}) &= V_j^+,  \quad &j \in \mathcal{B}, h \in \mathcal{H}, & \\
\sum_{t \in \mathcal{T}} (f_{i-1,i,t}^{l,h}+f_{i+1,i,t}^{l,h}) &= V_i^-,  \quad & i \in \mathcal{S}, h \in \mathcal{H},  \\
\sum_{\substack{ t \in \mathcal{T} }} (f_{\delta(k)-1,k,t}^{w,h} + f_{\delta(k)+1,k,t}^{w,h}) &= V_k^-,  \quad & k \in \mathcal{W}, h \in \mathcal{H}.
\end{align*}

The borrow and waste pit constraints ensure that the total sum of borrow/waste flows from/to a borrow/waste pit must not exceed the capacity of that borrow/waste, i.e.,
\begin{align*} 
\sum_{\substack{ t \in \mathcal{T} }} (f_{j,\vartheta(j)+1,t}^{b,h} + f_{j,\vartheta(j)-1,t}^{b,h}) &\le C_j^b,  \quad j \in \mathcal{B}, h \in \mathcal{H},  \\
\sum_{\substack{ t \in \mathcal{T} }} (f_{\delta(k)-1,k,t}^{w,h} + f_{\delta(k)+1,k,t}^{w,h}) &\le C_k^w,  \quad k \in \mathcal{W}, h \in \mathcal{H}.
\end{align*}

Before introducing block constraints we require some additional notation.  Let $\bar{\mathcal{I}}^2 =\{ (k_1,k_2) \in \mathcal{I} \times \mathcal{I}: k_1 < k_2, \text{ and } k_1, k_2$ are two consecutive blocks with no access road in between them\}, $\bar{\mathcal{I}}_\leftarrow =\{ k \in \mathcal{I}:$ the set of blocks $k$ before which there is no access road\}, and $\bar{\mathcal{I}}_\rightarrow =\{ k \in \mathcal{I}:$ the set of blocks $k$ after which there is no access road\}. The parameter $M_i$ is used to model binary constraints; in theory it can be any large enough number but for numerical reasons it is set to the largest possible cut or fill volume at a section $i$. For the sake of simplicity, we use $M$ instead of $M_i$.

Earth movement is not allowed between two sections when there is a block between them. This means that for each hauling path $h \in \mathcal{H}$, if a section $i$ is a block then material cannot be moved over that block or section. In other words, the whole material moved to the section $i$ from the section $i-1$ must be loaded to section $i$; i.e., for the left transit nodes $\forall h \in \mathcal{H},$ $f_{i-1,i,t}^{r,h} = f_{i-1,i,t}^{l,h}$. Moreover, if the section $i$ is a block the only material that can move from this section to the next section ($i+1$) is the material that has been excavated from section $i$, i.e., for the left transit nodes $\forall h \in \mathcal{H}, $ $f_{i,i+1,t}^{r,h} = f_{i,i+1,t}^{u,h}$. In the case of the right transit nodes, $f_{i+1,i,t}^{r,h} = f_{i+1,i,t}^{l,h}$ and $f_{i,i-1,t}^{r,h} = f_{i,i-1,t}^{u,h}$ ensure that no material can move across section $i$ to the previous section. Therefore, for the left transit nodes, the constraints
\begin{align} \label{14}
-(f_{\gamma (i)-1, \gamma (i) ,t}^{r,h} - f_{\gamma (i)-1,\gamma (i),t}^{l,h}) &\le My_{i,t-1}, \quad \forall i \in \mathcal{I}, h \in \mathcal{H}, t \in \mathcal{T} , \\
 \label{15}
f_{\gamma (i)-1, \gamma (i) ,t}^{r,h} - f_{\gamma (i)-1,\gamma (i),t}^{l,h} &\le My_{i,t-1}, \quad \forall i \in \mathcal{I}, h \in \mathcal{H}, t \in \mathcal{T} , \\
 \label{16}
-(f_{\gamma (i), \gamma (i)+1 ,t}^{r,h} - f_{\gamma (i),\gamma (i)+1,t}^{u,h}) &\le My_{i,t-1}, \quad \forall i \in \mathcal{I}, h \in \mathcal{H}, t \in \mathcal{T} ,\\
 \label{17}
f_{\gamma (i), \gamma (i)+1 ,t}^{r,h} - f_{\gamma (i),\gamma (i)+1,t}^{u,h}) &\le My_{i,t-1}, \quad \forall i \in \mathcal{I}, h \in \mathcal{H}, t \in \mathcal{T} ,
\end{align}
ensure that no material movement is allowed across a block. In the same way, for the right transit nodes the constraints
\begin{align} \label{18}
-(f_{\gamma (i)+1, \gamma (i) ,t}^{r,h} - f_{\gamma (i)+1,\gamma (i),t}^{l,h}) &\le My_{i,t-1}, \quad \forall i \in \mathcal{I}, h \in \mathcal{H}, t \in \mathcal{T} ,\\
 \label{19}
f_{\gamma (i)+1, \gamma (i) ,t}^{r,h} - f_{\gamma (i)+1,\gamma (i),t}^{l,h} &\le My_{i,t-1}, \quad \forall i \in \mathcal{I}, h \in \mathcal{H}, t \in \mathcal{T} ,\\
 \label{20}
-(f_{\gamma (i), \gamma (i)-1 ,t}^{r,h} - f_{\gamma (i),\gamma (i)-1,t}^{u,h}) &\le My_{i,t-1}, \quad \forall i \in \mathcal{I}, h \in \mathcal{H}, t \in \mathcal{T} ,\\
 \label{21}
f_{\gamma (i), \gamma (i)-1 ,t}^{r,h} - f_{\gamma (i),\gamma (i)-1,t}^{u,h} &\le My_{i,t-1}, \quad \forall i \in \mathcal{I}, h \in \mathcal{H}, t \in \mathcal{T} ,
\end{align}
ensure no material movement across a block.

Material movement is also not allowed between two blocks, before the first block, and after the last block with no access roads, until the blocks are removed. Thus, for all $h \in \mathcal{H}$, the transit, borrow, and waste flows are forbidden between two blocks with no access roads, until one of the blocks is removed, i.e.,
\begin{align*}
f_{i,i+1,t}^{r,h} &\le M(y_{k_1,t-1}+y_{k_2,t-1}), \\ &\quad \forall i \in \mathcal{S}, (k_1,k_2) \in \bar{\mathcal{I}}^2, h \in \mathcal{H}, t \in \mathcal{T}, \gamma(k_1) \le i, i+1 \le \gamma(k_2), \\
f_{i+1,i,t}^{r,h} &\le M(y_{k_1,t-1}+y_{k_2,t-1}), \\ &\quad \forall i \in \mathcal{S}, (k_1,k_2) \in \bar{\mathcal{I}}^2, h \in \mathcal{H}, t \in \mathcal{T}, \gamma(k_1) \le i, i+1 \le \gamma(k_2), \\
f_{j,\vartheta(j)+1,t}^{b,h} &\le M(y_{k_1,t-1}+y_{k_2,t-1}), \\ &\quad \forall j \in \mathcal{B}, (k_1,k_2) \in \bar{\mathcal{I}}^2, h \in \mathcal{H}, t \in \mathcal{T}, \gamma(k_1) \le \vartheta(j)-1, \vartheta(j)+1 \le \gamma(k_2), \\
f_{j,\vartheta(j)-1,t}^{b,h} &\le M(y_{k_1,t-1}+y_{k_2,t-1}), \\ &\quad \forall j \in \mathcal{B}, (k_1,k_2) \in \bar{\mathcal{I}}^2, h \in \mathcal{H}, t \in \mathcal{T}, \gamma(k_1) \le \vartheta(j)-1, \vartheta(j)+1 \le \gamma(k_2), \\
f_{\delta(j)-1,j,t}^{w,h} &\le M(y_{k_1,t-1}+y_{k_2,t-1}), \\ &\quad \forall j \in \mathcal{W}, (k_1,k_2) \in \bar{\mathcal{I}}^2, h \in \mathcal{H}, t \in \mathcal{T}, \gamma(k_1) \le \delta(j)-1, \delta(j)+1 \le \gamma(k_2), \\
f_{\delta(j)+1,j,t}^{w,h} &\le M(y_{k_1,t-1}+y_{k_2,t-1}), \\ &\quad \forall j \in \mathcal{W}, (k_1,k_2) \in \bar{\mathcal{I}}^2, h \in \mathcal{H}, t \in \mathcal{T}, \gamma(k_1) \le \delta(j)-1, \delta(j)+1 \le \gamma(k_2).
\end{align*}

After the last block with no access road, the transit, borrow, and waste flows are not allowed until the block is removed:
\begin{align*} 
f_{i,i+1,t}^{r,h} &\le My_{k,t-1},\quad \forall i \in \mathcal{S}, k \in \bar{\mathcal{I}}_\rightarrow, h \in \mathcal{H}, t \in \mathcal{T}, \gamma(k) \le i, i+1 \le n, \\
f_{i+1,i,t}^{r,h} &\le My_{k,t-1},\quad \forall i \in \mathcal{S}, k \in \bar{\mathcal{I}}_\rightarrow, h \in \mathcal{H}, t \in \mathcal{T}, \gamma(k) \le i, i+1 \le n, \\
f_{j,\vartheta(j)+1,t}^{b,h} &\le My_{k,t-1}, \quad \forall j \in \mathcal{B}, k \in \bar{\mathcal{I}}_\rightarrow, h \in \mathcal{H}, t \in \mathcal{T}, \gamma(k) \le \vartheta(j)-1, \vartheta(j)+1 \le n, \\
f_{j,\vartheta(j)-1,t}^{b,h} &\le My_{k,t-1}, \quad \forall j \in \mathcal{B}, k \in \bar{\mathcal{I}}_\rightarrow, h \in \mathcal{H}, t \in \mathcal{T}, \gamma(k) \le \vartheta(j)-1, \vartheta(j)+1 \le n,\\
f_{\delta(j)-1,j,t}^{w,h} &\le My_{k,t-1}, \quad \forall j \in \mathcal{W}, k \in \bar{\mathcal{I}}_\rightarrow, h \in \mathcal{H}, t \in \mathcal{T}, \gamma(k) \le \delta(j)-1, \delta(j)+1 \le n, \\
f_{\delta(j)+1,j,t}^{w,h} &\le My_{k,t-1}, \quad \forall j \in \mathcal{W}, k \in \bar{\mathcal{I}}_\rightarrow, h \in \mathcal{H}, t \in \mathcal{T}, \gamma(k) \le \delta(j)-1, \delta(j)+1 \le n.
\end{align*}

In the same way, before the first block with no access road, the transit, borrow, and waste flows are not allowed until the block is removed:
\begin{align*} 
f_{i,i+1,t}^{r,h} &\le My_{k,t-1},\quad \forall i \in \mathcal{S}, k \in \bar{\mathcal{I}}_\leftarrow, h \in \mathcal{H}, t \in \mathcal{T}, 1 \le i, i+1 \le \gamma(k), \\
f_{i+1,i,t}^{r,h} &\le My_{k,t-1},\quad \forall i \in \mathcal{S}, k \in \bar{\mathcal{I}}_\leftarrow, h \in \mathcal{H}, t \in \mathcal{T}, 1 \le i, i+1 \le \gamma(k), \\
f_{j,\vartheta(j)+1,t}^{b,h} &\le My_{k,t-1}, \quad \forall j \in \mathcal{B}, k \in \bar{\mathcal{I}}_\leftarrow, h \in \mathcal{H}, t \in \mathcal{T}, 1 \le \vartheta(j)-1, \vartheta(j)+1 \le \gamma(k), \\
f_{j,\vartheta(j)-1,t}^{b,h} &\le My_{k,t-1}, \quad \forall j \in \mathcal{B}, k \in \bar{\mathcal{I}}_\leftarrow, h \in \mathcal{H}, t \in \mathcal{T}, 1 \le \vartheta(j)-1, \vartheta(j)+1 \le \gamma(k),\\
f_{\delta(j)-1,j,t}^{w,h} &\le My_{k,t-1}, \quad \forall j \in \mathcal{W}, k \in \bar{\mathcal{I}}_\leftarrow, h \in \mathcal{H}, t \in \mathcal{T}, 1 \le \delta(j)-1, \delta(j)+1 \le \gamma(k), \\
f_{\delta(j)+1,j,t}^{w,h} &\le My_{k,t-1}, \quad \forall j \in \mathcal{W}, k \in \bar{\mathcal{I}}_\leftarrow, h \in \mathcal{H}, t \in \mathcal{T}, 1 \le \delta(j)-1, \delta(j)+1 \le \gamma(k).
\end{align*}

After excavating/embanking the required amount of earth from/to a section with a block, the block is considered removed. The \textit{block removal indicator constraints} are used to satisfy this expectation as follows
\begin{align*} 
\sum_{t =0}^{u} (f_{\gamma(k),\gamma(k)+1,t}^{u,h}+f_{\gamma(k),\gamma(k)-1,t}^{u,h})+M(1-y_{k,u}) &\ge V_{\gamma(k)}^+, \quad \forall k \in \mathcal{I}, h \in \mathcal{H}, u \in \mathcal{T},  \\
\sum_{t =0}^{u} (f_{\gamma(k)-1,\gamma(k),t}^{l,h}+f_{\gamma(k)+1,\gamma(k),t}^{l,h})+M(1-y_{k,u}) &\ge V_{\gamma(k)}^-, \quad \forall k \in \mathcal{I}, h \in \mathcal{H}, u \in \mathcal{T}.
\end{align*}

At the end of each time step $t \in \mathcal{T}$, at least one block should be removed so that removing all the blocks does not take more than $n_b+1$ time steps. The \textit{block removal enforcement} constraint guarantees this,
\[ 
\sum_{t =0}^{u}\sum_{k \in \mathcal{I}} y_{k,t} \ge t+1 \quad \forall u \in \setminus \{n_b\}.
\]
In addition, when a block is removed it should remain removed. This is done by the \textit{monotonicity constraint} as follows
\[ 
y_{k,t} \ge y_{k,t-1}, \quad \forall k \in \mathcal{I}, t \in \setminus \{0\}.
\]

The continuity constraints,
\begin{align*}
 P_{g-1}(s_{\varphi(g,1)}) &= P_{g}(s_{\varphi(g,1)}) \quad \forall g \in \mathcal{G} \setminus \{1\}, \\
 P'_{g-1}(s_{\varphi(g,1)}) &= P'_{g}(s_{\varphi(g,1)}) \quad \forall g \in \mathcal{G} \setminus \{1\},
\end{align*}
ensure that the height and the slope of the first section of a segment is equal to the height and the slope of the last section of the previous segment.

The \textit{volume constraints}
\[ V_i^+ - V_i^- = A_i u_i, \]
where $A_i$ is the area of section $i$, guarantee that the total excavated/embanked volume of material from/to a section equals to the volume difference between the road profile and the ground profile for that section.

There are restrictions on the grade of the road profile to satisfy safety considerations. \textit{Slope constraints}
\[G_L \leq P'_{g}(s_{\varphi(g,1)}) \le G_U  \quad \forall g \in \mathcal{G} \setminus \{1\}, \]
are used to constrain the spline segments within a closed interval $[G_L,G_U]$, in which $G_L$ and $G_U$ are the minimum and the maximum valid grades, respectively.

Finally, We need \textit{bound constraints} to restrict the domain of variables with
\begin{align*} 
 0 \le V_i^+ &\le M, \quad \forall i \in \mathcal{S \cup B}, \\
 0 \le V_i^- &\le M, \quad \forall i \in \mathcal{S \cup W},
\end{align*}
and
\begin{align*}
 f_{i,i+1,t}^{r,h} &\ge 0, \quad f_{i,i-1,t}^{r,h} \ge 0 &\quad \forall i \in \mathcal{S}, h \in \mathcal{H}, t \in \mathcal{T}, \\
 f_{i-1,i,t}^{l,h} &\ge 0, \quad f_{i+1,i,t}^{l,h} \ge 0 &\quad \forall i \in \mathcal{S}, h \in \mathcal{H}, t \in \mathcal{T}, \\
 f_{i,i+1,t}^{u,h} &\ge 0, \quad f_{i,i-1,t}^{u,h} \ge 0 &\quad \forall i \in \mathcal{S}, h \in \mathcal{H}, t \in \mathcal{T}, \\
 f_{j,\vartheta(j)+1,t}^{b,h} &\ge 0, \quad f_{j,\vartheta(j)-1,t}^{b,h} \ge 0 &\quad \forall j \in \mathcal{B}, h \in \mathcal{H}, t \in \mathcal{T}, \\
 f_{\delta(j)-1,j,t}^{w,h} &\ge 0, \quad f_{\delta(j)+1,j,t}^{w,h} \ge 0 &\quad \forall j \in \mathcal{W}, h \in \mathcal{H}, t \in \mathcal{T}.
\end{align*}

\subsection{The CTG and QNF Models}

To evaluate the performance of the MH-QNF model, we compare the efficiency, robustness, and accuracy of the proposed model with the CTG and QNF models presented in \cite{4}.  At this point, the QNF model of \cite{4} is quite easy to describe -- the QNF model is a special case of the MH-QNF model, where the number of haul types is one ($n_h =1$).  The CTG model, which appears in \cite{4, 5}, is slightly harder to outline.  In essence, instead of hauling paths, for every pair of nodes $(i, j)$, the CTG model places an arc moving from $i$ to $j$.  The flow through the arc is a continuous variable and the arc cost is computed based on $i$ and $j$.  The CTG model allows a great deal of flexibility in arc-costs, but the resulting number of variables is very large.  In particular, in the MH-QNF model (and hence the QNF model) the total number of variables grows linearly with the number of sections.  Conversely, for the CTG model, the total number of variables grows quadratically with the number of sections.   As the number of sections is typically the dominating factor in these models, the CTG model is typically much larger than the MH-QNF model.

\section{Numerical results}

To evaluate the performance of the MH-QNF model, we compare the efficiency, robustness and accuracy of the proposed model with both CTG and QNF models presented in \cite{4} in terms of timing and optimal cost.

\begin{table}[ht]
\caption{The specifications of the basic roads used to assemble the test collection.}\label{Table0}
\begin{center}
\begin{tabularx} {0.85\textwidth}{l c c c}
  \toprule 
  \head{Road} & \head{Length (km)} & \head{Section Length (m)} & \head{Number of sections}  \\
  \toprule 
  A & 1 & 20 & 50 \\
  B & 5 & 100 & 50 \\
  C & 2 & 20 & 100 \\
  D & 3 & 20 & 150 \\
  E & 15 & 100 & 150 \\
  F & 20 & 100 & 200 \\
  G & 9 & 20 & 450 \\
  \bottomrule 
\end{tabularx}

\end{center}
\end{table}

The selected problem test collection consists of $60$ problems. These problems are generated by changing various parameters (i.e., number of sections, blocks, access roads, offset levels and section lengths) in 7 distinct road samples, provided by our industrial partner (see Table \ref{Table0}).  All tests represent realistic road design problems. 

We consider 5 hours timeout for all the experiments; beyond this a test is reported as unsuccessful. Linear programming feasibility tolerance and relative mixed integer programming gap tolerance parameters are set to $10^{-06}$ and 1\%, respectively. The upper bound and lower bound for slope constraints are 0.1 and -0.1 ($\pm 10\%$). All other parameters are set to their default value.

After selection of the test collection, we solved the problems by different models using the same computational environment. The workstation used for experiments has an Intel(R)Xeon(R) CPU W3565 3.20GHz (4 cores and 8 threads) processor and 24.0 GB of RAM. To solve the developed mathematical models we use the academic version of the IBM ILOG CPLEX optimizer 12.51 edition \url{http://www.cplex.com}.

One of the influencing factors on our experiments is the cost components used in the MH-QNF model shown in Table \ref{Table1} (specific values were selected in consultation with our industry partner). Excavation and embankment cost components depend on the different types of materials. Loading and hauling cost components depend on the type of hauling path. The Avg Value column is the average of the cost values for a cost component. All the three models are multi material and share the same excavation and embankment cost components for each material. Both MH-QNF model and the CTG model use multiple hauling paths with different hauling and loading costs. Therefore, in the case of the QNF model which has a single hauling path, for the loading and hauling cost components we need to test the model with different pairs of hauling and loading cost components.

\begin{table}
\captionof{table}{Cost components in $\$ \slash cu.m$.\label{Table1}}
\begin{center}

\begin{tabularx} {0.7\textwidth}{l c c c}
  \toprule[1.5pt]
  \head{Cost component} & \head{Cost Type} & \head{Value} &   \\
  \toprule[1.5pt]
  Excavation cost & M1 & 4.000 &  \\
  				  & M2 & 4.000 &  \\
   				  & M3 & 20.000 &  \\
  				  & M4 & 4.000 &  \\   \midrule
  				
  Embankment cost & M1 & 2.000 &  \\
    			  & M2 & 2.000 &  \\
     			  & M3 & 1.800 &  \\
    			  & M4 & 2.000 &  \\
  \toprule[1.5pt]
  \head{Cost component} & \head{Cost Type} & \head{Value} & \head{Avg Value}  \\
  \toprule[1.5pt]    			
  Loading cost    & Short & 0.000 & 1.067 \\
    			  & Middle & 0.600 &  \\     			
    			  & Long & 2.600 &  \\ \midrule
    			
  Hauling cost    & Short & 0.008 & 0.005 \\
    			  & Middle & 0.004 &  \\     			
    			  & Long & 0.002 &  \\
  \bottomrule[1.5pt]
\end{tabularx}

\end{center}
\end{table}

Since the MH-QNF model is an extension of the QNF model, to get a better understanding of the importance of loading and hauling cost components, in addition to the average values of loading and hauling costs (1.067, 0.005, respectively), we test the QNF model with the other three pairs of cost values for loading and hauling cost components reported for different hauls (Table \ref{Table1}); i.e., we run the QNF model four times with different loading and hauling costs: average case (1.067, 0.005), short haul (0.000, 0.008), middle haul (0.600, 0.004), and long haul (2.600, 0.002). Thus, we have 4 different setups for the QNF model as follows.

\begin{itemize}
\item  \textbf{QNF-S}: the model QNF configured with short hauling and loading cost values.
\item \textbf{QNF-M}: the model QNF configured with middle hauling and loading cost values.
\item \textbf{QNF-L}: the model QNF configured with long hauling and loading cost values.
\item \textbf{QNF-A}: the model QNF configured with average hauling and loading cost values as shown in Table~\ref{Table1}.
\end{itemize}

\begin{figure} [ht]
    \centering
	\includegraphics[scale=0.4]{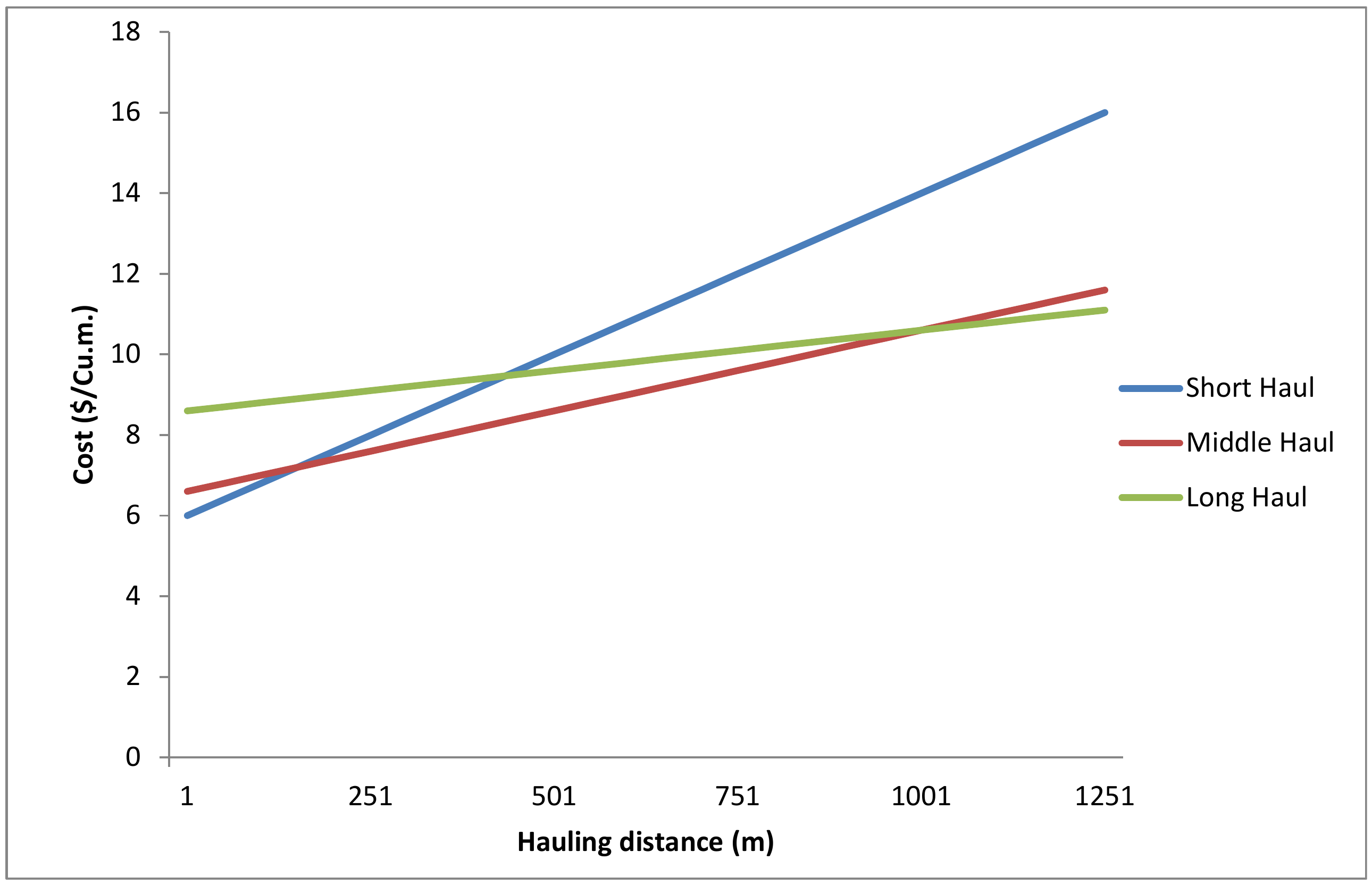}
    \caption{Cost configuration.}
    \label{7}
\end{figure}

As shown in Figure \ref{7}, to evaluate the efficiency of using a variety of hauling paths over a single path in a road construction site, the loading and hauling cost components are selected such that
\begin{itemize}
\item if \textbf{$\text{ hauling distance} < 150 m$} then short hauling path will be picked as the optimum method of transferring material;
\item if \textbf{$ 150m \le \text{hauling distance} < 1000 m$} then middle hauling path will be picked;
\item and finally, if \textbf{$\text{hauling distance} \ge 1000 m$} then the long hauling path will be selected.
\end{itemize}
(These specific values were selected after discussion with our industrial partner.)

In \cite{4} the authors propose six potential techniques to improve the speed of solution search. They extensively test all the techniques with both models and identify the setup which works the best for both models. Their results highlight two techniques as the most efficient methods. Therefore, we test all models using two different configurations. In the first configuration, we use SOS2 variables for the volume constraints and SOS1 variables for modeling the block constraints in all the models. In the second configuration, again we use SOS2 variables for the volume constraints but for the block constraints the basic technique is employed.

In general, in our experiments we use 12 model configurations. To increase the readability of the report on the results we use the following naming convention.
\begin{enumerate}
\item Models that use the \textit{basic} technique for blocks:
\begin{itemize}
\item \textbf{MQN-B}: The MH-QNF model.
\item \textbf{CTG-B}: The CTG model with multi haul cost components.
\item \textbf{QNS-B}: The QNF-S model.
\item \textbf{QNM-B}: The QNF-M model.
\item \textbf{QNL-B}: The QNF-L model.
\item \textbf{QNA-B}: The QNF-A model.
\end{itemize}
\item Models that use the \textit{SOS1} techniques for blocks:
\begin{itemize}
\item \textbf{MQN-S1}: The MH-QNF model.
\item \textbf{CTG-S1}: The CTG model with multi haul cost components.
\item \textbf{QNS-S1}: The QNF-S model.
\item \textbf{QNM-S1}: The QNF-M model.
\item \textbf{QNL-S1}: The QNF-L model.
\item \textbf{QNA-S1}: The QNF-A model.
\end{itemize}
\end{enumerate}

\subsection{Reporting the results}

As a first pass at examining model accuracy, we summarize the relative errors in cost values obtained by all models in Table \ref{tab:errors}. The column ``opt. found'' is the number of problems where the model optimum value was found before timeout occurred (regardless of relative error). The columns ``min./mean/max. error'' are the minimum/mean/maximum relative error over all problems solved when compared to CTG configurations (CTG-B or CTG-S1), based on only problems solved by both CTG configuration and the model examined.

\begin{table}[ht]
\begin{center}
\caption{Summary of model accuracy for each model.}
\label{tab:errors}
\begin{tabular}{c l c c c c}\hline

\toprule 
\head{\#} & \head{Model} & \head{opt. found} & \head{min. error (\%)} & \head{mean error (\%)}& \head{max. error (\%)}  \\
\toprule 
1 & CTG-B & 51 & -- & -- & -- \\
2 & MQN-B & 56 & 0.00 & 0.47 & 2.20 \\
3 & QNS-B & 56 & 27.15 & 57.10 & 100.92 \\
4 & QNM-B & 56 & 12.33 & 25.65 & 40.74 \\
5 & QNL-B & 56 & 18.93 & 25.81 & 37.95 \\
6 & QNA-B & 56 & 22.19 & 40.22 & 62.98 \\ \midrule
7 & CTG-S1 & 45 & -- & -- & -- \\
8 & MQN-S1 & 57 & 0.00 & 0.37 & 1.43 \\
9 & QNS-S1 & 56 & 27.15 & 57.56 & 100.92 \\
10 & QNM-S1 & 56 & 12.45 & 25.78 & 40.74 \\
11 & QNL-S1 & 56 & 18.93 & 25.81 & 37.39 \\
12 & QNA-S1 & 56 & 22.36 & 40.44 & 62.98 \\
\toprule 
\end{tabular}\end{center}
\end{table}

Among the above 12 model configurations, none of the QNF model configurations are able to solve a single problem satisfactorily. We may therefore discard these methods as inaccurate.

Table \ref{Table2} shows the solution times for the remaining 4 configurations applied to 60 test problems when only 1\% relative error is allowed. ``NaN'' values indicate that the solver exceeded the 5 hours timeout.

\begin{center}
\begin{longtable}{c r r r r}
\caption{Solution times (in second) for 1\% relative error.} \label{Table2} \\

  \head{Test \#} & \head{MQN-B} & \head{CTG-B} & \head{MQN-S1} & \head{CTG-S1} \\
  \toprule[1.5pt]
1 & 23.84 & 19.02 & 15.21 & 33.71\\
2 & 11.00 & 27.18 & 15.28 & 35.41\\
3 & 26.51 & 190.64 & 20.51 & 237.69\\
4 & 168.62 & 7214.67 & 172.48 & 4848.31\\
5 & 42.07 & 434.39 & 52.39 & 452.03\\
6 & 13.05 & 232.84 & 28.26 & 1014.80\\
7 & 874.97 & 603.49 & 581.46 & 1928.96\\
8 & NaN & 12605.20 & 2409.50 & NaN\\
9 & 19.93 & 169.51 & 78.17 & 309.46\\
10 & 13.40 & 297.75 & 17.13 & 141.41\\
11 & 66.34 & 2348.93 & 60.28 & NaN\\
12 & 146.36 & 3929.18 & 124.01 & 9853.74\\
13 & 4.96 & 86.70 & 11.38 & 33.63\\
14 & 2.15 & 39.36 & 4.97 & 17.24\\
15 & 18.75 & 392.67 & 35.29 & 451.63\\
16 & 319.90 & NaN & 366.57 & NaN\\
17 & 13.81 & 81.61 & 23.76 & 174.21\\
18 & 238.37 & 1369.81 & 383.95 & 4559.33\\
19 & 30.17 & 799.55 & 68.35 & 1657.61\\
20 & 81.99 & 632.69 & 1175.49 & 4146.94\\
21 & 455.83 & NaN & 2743.00 & NaN\\
22 & 223.27 & 8901.20 & 1426.83 & NaN\\
23 & 65.00 & 1306.09 & 553.65 & 2437.76\\
24 & 3.65 & 91.19 & 7.87 & 55.34\\
25 & 9.64 & 750.86 & 25.49 & 314.50\\
26 & 76.04 & 5486.19 & 158.04 & NaN\\
27 & 6.87 & 353.50 & 14.64 & 57.91\\
28 & 6.16 & 493.60 & 14.19 & 248.43\\
29 & 5.52 & 67.06 & 9.20 & 19.16\\
30 & 30.44 & 2454.93 & 51.14 & 2432.61\\
31 & 9.89 & 179.33 & 9.83 & 35.04\\
32 & 18.15 & 619.50 & 78.73 & 241.43\\
33 & 11.03 & 1361.16 & 25.27 & 675.14\\
34 & 30.08 & 2385.43 & 91.46 & NaN\\
35 & 22.71 & 2463.64 & 35.74 & 2543.56\\
36 & 120.08 & NaN & 273.13 & 9350.41\\
37 & 57.23 & 1702.74 & 90.71 & NaN\\
38 & 3.12 & 4.96 & 4.95 & 4.98\\
39 & 17.09 & 230.50 & 36.15 & 623.73\\
40 & 18020.00 & NaN & NaN & NaN\\
41 & 3361.52 & NaN & 9224.50 & NaN\\
42 & 392.59 & NaN & 455.71 & NaN\\
43 & 308.46 & NaN & 313.39 & NaN\\
44 & NaN & NaN & NaN & NaN\\
45 & 70.47 & 2014.47 & 168.04 & NaN\\
46 & 50.75 & 1289.76 & 122.84 & 3528.95\\
47 & 8.85 & 7.44 & 17.43 & 7.78\\
48 & 5.70 & 9.02 & 10.61 & 9.76\\
49 & 3.22 & 1.87 & 6.72 & 1.89\\
50 & 1.36 & 0.94 & 2.46 & 0.99\\
51 & 4.83 & 6.11 & 7.72 & 6.14\\
52 & 8.97 & 28.74 & 14.42 & 28.87\\
53 & 14.37 & 40.89 & 23.76 & 41.09\\
54 & 13.58 & 27.02 & 21.90 & 27.86\\
55 & 13.68 & 30.44 & 23.01 & 31.26\\
56 & 3.78 & 4.74 & 6.34 & 4.68\\
57 & 73.50 & 1373.68 & 94.70 & 1361.67\\
58 & NaN & NaN & NaN & NaN\\
59 & 15.24 & 122.75 & 15.90 & 128.76\\
60 & 6.99 & 14.09 & 8.58 & 16.29\\
\bottomrule[1.5pt]
\end{longtable}
\end{center}

The reliability and robustness of an optimization method is defined as the ability of the method to perform well over a wide range of optimization problems. To evaluate the reliability of the model we use two measures: success rate and computational accuracy. The success rate is defined as the number of problems in the given problem set that are successfully solved to optimality by the optimization method. We use performance profiles \cite{36} to evaluate the reliability and robustness of the models. While performance profiles are now prevalent in optimization benchmarking, we provide a brief review for completeness.  The performance profile of a solver $s$ is defined as
\[
\rho_s({\alpha})=\frac{1}{\lvert \mathcal{P} \rvert} \text {size} \{\mathit{p}  \in \mathcal{P}: r_{p,s} \le \alpha \},
\]
where $\lvert \mathcal{P} \rvert$ represents the cardinality of the problem set $\mathcal{P}$ and $r_{p,s}$ is the performance ratio defined as
\[
  r_{p,s}=
  \begin{cases}
    \frac{t_{p,s}}{\min\{t_{p,s}:s \in \mathcal{S}\}} &\text{if  the convergence test passed } \\
    \infty &\text{if the convergence test failed}
  \end{cases}
\]
in which $t_{p,s}$ is the performance measure which is the computational time in our case. In this equation, for a specific problem and the best solver in terms of $t_{p,s}$ has $r_{p,s} = 1$.

One of the advantages of performance profiles, is that they implicitly includes the performance ratio \textit{success rate} as a reliability factor. The value of $\rho_s({\alpha})$ gives a sense of how promising the solutions found by an optimization algorithm are relative to the best solution found by all the optimization algorithms that are compared together.

For each solution obtained for problem number $s$ using the model $m$ the percent relative error is calculated using the following formula:
\begin{equation} \label{53}
 E_{s,m}= \dfrac{Obj_s - Obj_b}{Obj_b},
\end{equation}
in which $Obj_s$ is the optimal objective value obtained for the problem $s$ by the model $m$ and $Obj_b$ is the objective value obtained by the benchmark model $b$. In our experiments, a solution found by model $m$ is considered \textit{successful} if $\lvert E_{s,m}\rvert \le 1.0\% $.

\begin{figure}[H]
    \centering
	\includegraphics[scale=0.50]{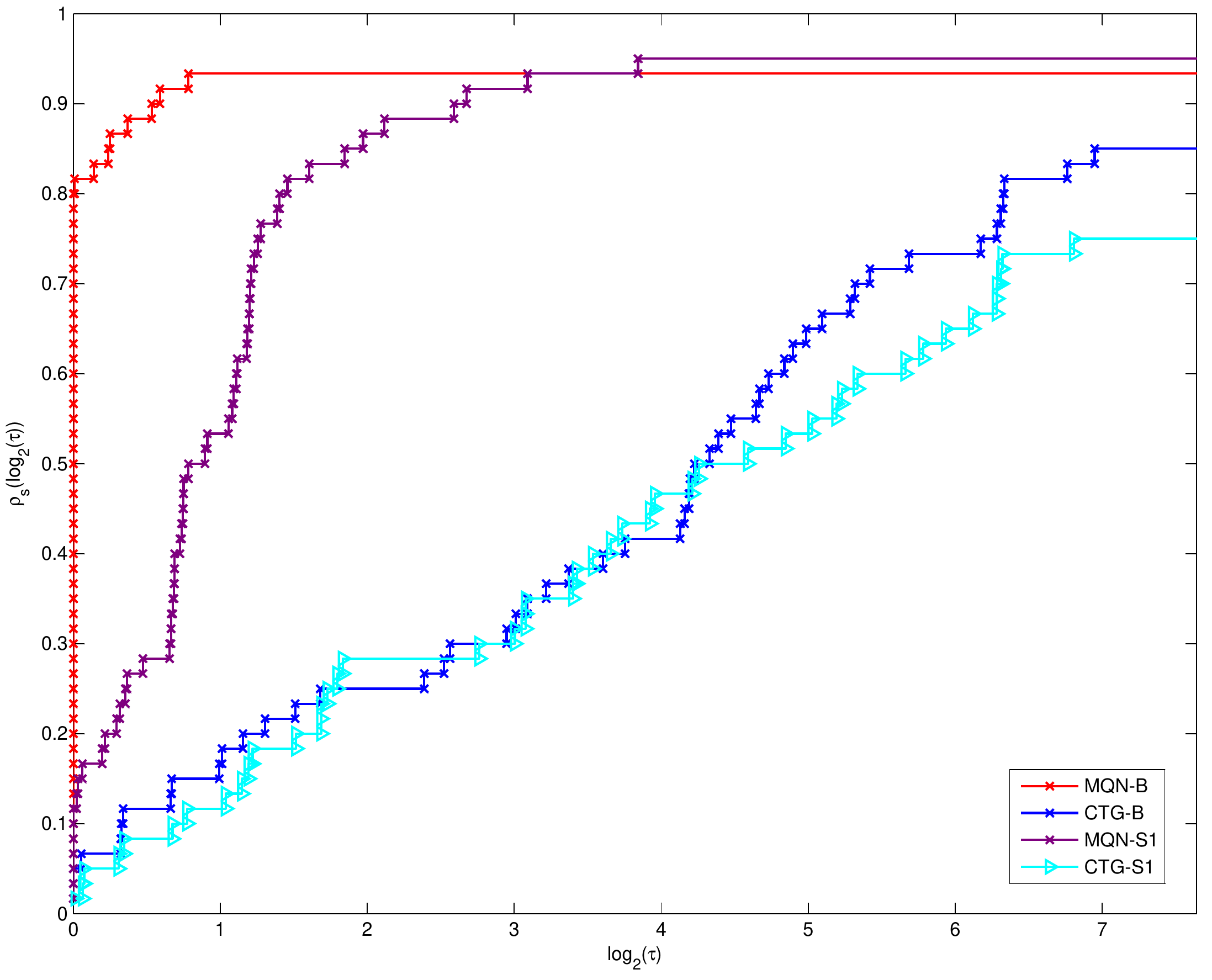}
    \caption{Performance profiles of the models with up to $1\%$ relative error accepted.}
    \label{5}
\end{figure}

When available, we consider the CTG model as the benchmark and we calculate the relative error for the solutions obtained using the MH-QNF model.  When the CTG model does not provide a solution but the MH-QNF does, we consider the MH-QNF relative error to be $0.0\%$.  Figure \ref{5} provides the resulting performance profile.

Examining Figure \ref{5}, it is clear that MH-QNF greatly outperforms CTG.  In terms of the block technique, the CTG model has a relatively better performance when it uses the basic technique for blocks by solving almost $82\%$ of the problems compared to the SOS1 version (CTG-S1) which ends up with solving $75 \%$. Similarly, in the case of the MH-QNF model, the performance of the model using the basic technique (MQN-B) is considerably higher than the SOS1 technique (MQN-S1). Except in one test, MQN-B solves most of the problems faster than MQN-S1. In general, among all the configurations, MQN-B is the winner and has a better performance than the CTG configurations.

\section{Conclusion}

In this research, we extend a model based on a previous quasi network flow model by \cite{4} for the vertical alignment problem considering the earthwork cost. In particular, we add multiple hauling flows to the model in order to improve the accuracy of the model and provide more modeling flexibility for the users. The model considers features such as blocks and side-slopes. We compare the model with two state-of-the-art models in the field: CTG and QNF models.

The CTG model allows a great deal of flexibility in arc-costs, but the resulting number of variables is very large. The total number of variables grows quadratically with the number of sections. As the number of sections is typically the dominating
factor in these models, the CTG model is typically much larger than the QNF and MH-QNF models.

In the QNF model, the total number of variables grows linearly with the number of sections.  However, in QNF model there is a linear relation between the hauling cost and the distance traveled by the volume of materials (e.g., the cost of moving 1 unit of earth for 100m equals to 100 times the cost of moving 1 unit for 1m), which is not true in real world. Therefore, we proposed the MH-QNF model to integrate the low computational cost of QNF model with the accuracy of CTG in solving the vertical alignment problem. The experimental results confirm that the new model is both highly accurate and considerably faster than the CTG model.

These results also provide an important advancement to solving the complete road design problem (where earthwork, vertical alignment, and horizontal alignment are optimized simultaneously). Recent research by \cite{MONDAL-15}, has viewed the complete road design problem as a bi-level optimization problem, where the horizontal alignment is optimized using a non-gradient based method which calls the vertical alignment problem as a black-box.  For this approach to work, it is critical that the vertical alignment is solved both quickly and accurately.  The MH-QNF model provides both of these attributes.

\section*{Acknowledgments}
This work was supported by Natural Sciences and Engineering Research Council of Canada (NSERC) under Collaborative Research and Development (CRD) Grant \#CRDPJ 411318-10. The Grant was sponsored by Softree Technical Systems Inc. This work was supported by Discovery Grants \#355571-2013 (Hare) and \#298145-2013 (Lucet) from NSERC. Part of the research was performed in the Computer-Aided Convex Analysis (CA2) laboratory funded by a Leaders Opportunity Fund (LOF) from the Canada Foundation for Innovation (CFI) and by a British Columbia Knowledge Development Fund (BCKDF).

The authors acknowledge the original ideas and concepts of Dr. Donovan Hare that contributed to the quasi network flow model.

\bibliographystyle{te} 

\bibliography{MHC-QNFModel}

\end{document}